\documentclass[pdflatex,sn-mathphys-num]{sn-jnl}

\usepackage{graphicx}%
\usepackage{multirow}%
\usepackage{amsmath,amssymb,amsfonts}%
\usepackage{amsthm}%
\usepackage{mathrsfs}%
\usepackage[title]{appendix}%
\usepackage{xcolor}%
\usepackage{textcomp}%
\usepackage{manyfoot}%
\usepackage{booktabs}%
\usepackage{algorithm}%
\usepackage{algorithmicx}%
\usepackage{algpseudocode}%
\usepackage{listings}%

\theoremstyle{thmstyleone}%
\newtheorem{theorem}{Theorem}
\newtheorem{proposition}[theorem]{Proposition}%

\theoremstyle{thmstyletwo}%
\newtheorem{example}{Example}%
\newtheorem{remark}{Remark}%

\theoremstyle{thmstylethree}%

\newtheorem{corollary}[theorem]{Corollary}
\newtheorem{lemma}[theorem]{Lemma} 
\newtheorem{problem}[theorem]{Problem}

\begin{document}

\title[On weak${^*}$-basic sequences]{On weak${^*}$-basic sequences in duals and biduals of spaces $C(X)$ and Quojections\vspace{0.5cm}

\small{\textit{To the memory of our friend Professor Steven A. Saxon.}}}


\author[1]{\fnm{Jerzy} \sur{K{\c{a}}kol}}\email{kakol@amu.edu.pl}\equalcont{These authors contributed equally to this work.}

\author*[2]{\fnm{Manuel} \sur{Lopez-Pellicer}}\email{mlopezpe@mat.upv.es}
\equalcont{These authors contributed equally to this work.}

\author[3]{\fnm{Wies{\l}aw} \sur{\'Sliwa}}\email{sliwa@amu.edu.pl; wsliwa@ur.edu.pl}
\equalcont{These authors contributed equally to this work.}

\affil*[1]{\orgdiv{Faculty of Mathematics and Informatics.}, \orgname{A. Mickiewicz University}, \orgaddress{\street{Uniwersytetu Poznańskiego 4}, \city{Pozna\'{n}}, \postcode{61-614}, \country{Poland}}}

\affil[2]{\orgdiv{Institute of Pure and Applied Mathematics}, \orgname{Universitat Politècnica de València}, \orgaddress{\street{Camino de Vera s/n}, \city{Valencia}, \postcode{46022}, \country{Spain}}}

\affil[3]{\orgdiv{Faculty of Exact and Technical Sciences}, \orgname{University of Rzesz\'ow}, \orgaddress{\street{al. Tadeusza Rejtana 16C}, \city{Rzesz\'ow}, \postcode{35-310}, \country{Poland}}}


\abstract{We show that for infinite Tychonoff spaces $X$ and $Y$ the weak$^\ast$-dual of $C_k(X\times Y)$ contains a basic sequence; moreover, the weak$^\ast$-bidual of $C_k(X)$ contains such a sequence as well.  When $X$ and $Y$ are infinite compact spaces, we single out a concrete sequence $(\mu_n)$ of finitely supported signed measures on $X\times Y$ with quantitative small-rectangle estimates, and we prove that every subsequence of $(\mu_n)$ admits a further subsequence which is strongly normal and forms a weak$^\ast$-basic sequence in the dual $C(X\times Y)^\ast$ of the Banach space $C(X\times Y)$. We also study the weak$^\ast$-basic sequence problem for Fr\'echet locally convex spaces in the class of quojections, and prove that for every quojection $E$ the bidual $E^{\ast\ast}$ admits a weak$^\ast$-basic sequence, while a long-standing open problem asks whether the dual of every infinite-dimensional Banach space admits a basic sequence in the weak$^{*}$-topology. Several examples and open questions are included, in particular for spaces $C(X)$ and for inductive limits of Fr\'echet spaces.}




\maketitle

\section{Introduction}\label{sec1}

For a locally convex space (lcs) $E$ we set $E^\ast$ for its dual and
use $\sigma(E,E^\ast)$ and $w^{*}=\sigma(E^\ast,E)$ to denote the weak and weak$^\ast$ topologies, respectively.  The strong dual of $E$ we denote by $E_b^\ast$, that is, $E^\ast$ endowed with the strong topology $\beta(E^\ast,E)$.

A sequence $(e_n)$ in $E$ is a \emph{Schauder basis} if every $x\in E$ admits a unique representation $x=\sum_{n=1}^\infty b_n e_n$ with convergence in the given topology of $E$ and the coefficient functionals $b_n^\ast: E \to \mathbb{R}, b_n^\ast (x)=b_n$ are continuous for any $n\in \mathbb{N}$. A sequence is a \emph{basic sequence} if it is a Schauder basis of its closed linear span.

The existence of weak$^\ast$-basic sequences in duals is naturally tied to the separable-quotient problem. For Fr\'echet spaces (i.e. metrizable and complete lcs), the classical theorem of Eidelheit implies that every non-normable Fr\'echet space has a quotient isomorphic to $\mathbb{R}^{\mathbb{N}}$; hence the dual of such a space contains a weak$^\ast$-basic sequence.  In contrast, the corresponding statement for Banach spaces remains unsettled: it is open whether the dual of every infinite-dimensional Banach space necessarily contains a weak$^\ast$-basic sequence, and there are barrelled locally convex spaces (not of the form $C_k(X)$ or $C_p(X)$) for which the dual admits no weak$^\ast$-basic sequence, see references below.

In this paper, we establish weak$^\ast$-basic sequences in duals of several natural function spaces.  For any infinite Tychonoff spaces $X$ and $Y$ we prove that $C_k(X\times Y)^\ast$ contains a weak$^\ast$-basic sequence, and that $C_k(X)^{\ast\ast}$ contains one as well.  In the compact case, we exploit an explicit sequence of finitely supported signed measures $(\mu_n)$ on $X\times Y$ (obtained in \cite[Theorem~1]{Kakol-Sliwa}) satisfying quantitative small-rectangle estimates; this yields that every subsequence of $(\mu_n)$ admits a strongly normal subsequence which is weak$^\ast$-basic.

We also treat the class of quojections.  Recall that a Fr\'echet space $E$ is a \emph{quojection} if for every continuous seminorm $p$ on $E$ the quotient $E/\ker p$, endowed with the quotient topology, is a Banach space.  We show that for every quojection $E$ the bidual $E^{\ast\ast}$ contains a weak$^\ast$-basic sequence.

Finally, motivated by these results, we discuss inductive limits of Fr\'echet spaces.  If $E$ is the strict inductive limit of a sequence of Fr\'echet spaces, then $E^\ast$ (with the weak$^\ast$ topology) contains a complemented copy of $\mathbb{R}^{\mathbb{N}}$ and hence a weak$^\ast$-basic sequence.  This leads to the question whether the same conclusion holds for arbitrary $(LF)$-spaces.

\section{Preliminaries}

A long-standing problem asks whether the dual of every infinite-dimensional Banach space admits a basic sequence in the $w^{*}$-topology (see, for example, \cite{Sevilla} and the references, and  \cite{Mujica}). In the Fr\'echet setting this question is tightly linked to separable quotients and admits a particularly transparent characterisation.  Indeed, in \cite[Proposition 1]{S-W} \'Sliwa and W\'ojtowicz proved the following theorem, see also \cite{K-S}.
\begin{theorem}[\'Sliwa--W\'ojtowicz]\label{I}
For a lcs $E$ consider the following statements:
\begin{enumerate}
\item $E$ admits an infinite-dimensional separable quotient.
\item $E$ admits an infinite-dimensional separable quotient having a Schauder basis.
\item $(E^{*},\sigma(E^{*},E))$ has a basic sequence.
\item There exists in $E$ a strictly increasing sequence of closed vector subspaces whose union is dense in $E$.
\end{enumerate}
Then $(2)\Rightarrow(3)\Rightarrow(4)\Leftrightarrow(1)$. Moreover, all statements are equivalent if $E$ is a Fr\'echet lcs.
\end{theorem}
Theorem~\ref{I} extends the classical result of Johnson and Rosenthal \cite[Theorem III.1]{JR} for separable Banach spaces.
It also applies to every non-normable Fr\'echet lcs, because such a space has a quotient isomorphic to $\mathbb{R}^{\mathbb{N}}$ by \cite{Eidelheit}.
For Banach spaces the situation is substantially more delicate; nevertheless, Theorem~\ref{I} applies to various concrete classes, for instance, to the Lipschitz free space setting, cf. \cite{Hajek}.

\bigskip
\noindent\textbf{(A) Spaces of continuous functions.}
Our first main result identifies broad families of $C$-spaces whose duals (or biduals) contain weak$^\ast$-basic sequences.
\begin{theorem}\label{VVI}
Let $X$ and $Y$ be infinite Tychonoff spaces. Then:
\begin{enumerate}
\item The dual of $C_{k}(X\times Y)$ admits a weak$^\ast$-basic sequence.
\item The bidual of $C_{k}(X)$ admits a weak$^\ast$-basic sequence.
\item If $C_{k}(X)$ is barrelled, then the dual of $C_{k}(X)$ admits a weak$^\ast$-basic sequence.
\item If $X$ contains an infinite discrete subset which is $C^{*}$-embedded in $X$, then the dual of $C_{k}(X)$ admits a weak$^\ast$-basic sequence.
\item If $X$ and $Y$ are compact and $(\mu_n)$ is a sequence as in Theorem~\ref{SU} below, then every subsequence $(\mu_{k_n})$ of $(\mu_n)$ contains a strongly normal subsequence that is a weak$^\ast$-basic sequence in $C(X\times Y)^{*}$.
\end{enumerate}
\end{theorem}
Item~(5) relies on a concrete sequence of finitely supported measures from (\cite[Theorem 1]{Kakol-Sliwa}); for convenience, we record the required statement here:
\begin{theorem}[K\c akol-\'Sliwa]\label{SU}
Let $X$ and $Y$ be infinite compact spaces and let $n\in\mathbb{N}$.
Let $K_n\times L_n$ be a finite subset of $X\times Y$ such that $|K_n|=2^{n}$ and $|L_n|=n$.
Let $\varphi_n:K_n\to\{-1,1\}^{L_n}$ be a bijection.
Then
\[
\mu_n:=\frac{1}{n2^{n}}\sum_{(s,j)\in K_n\times L_n} \varphi_n(s)(j)\,\delta_{(s,j)}
\]
is a finitely supported signed measure on $X\times Y$ such that
\[
(1)\;\; \|\mu_n\|=1 \ \text{and}\ \operatorname{supp}(\mu_n)=K_n\times L_n;
\]
\[
(2)\;\; \frac{1}{2\sqrt{\pi}}\frac{1}{\sqrt{n}}
<\sup_{A\times B\subset X\times Y} |\mu_n(A\times B)|
<\frac{2}{\sqrt{\pi}}\frac{1}{\sqrt{n}};
\]
\[
(3)\;\; |\mu_n(f\otimes g)|\leq \frac{8}{\sqrt{\pi}}\frac{1}{\sqrt{n}}\|f\otimes g\|_{\infty}, \mbox{for all}\;(f,g)\in C(X)\times C(Y),
\]
where $f\otimes g:X\times Y\to\mathbb{R}, (f\otimes g)(x,y)=f(x)g(y)$.

In particular, $\mu_n(h)\to_n 0$ for every $h\in C(X\times Y)$.
\end{theorem}
In Section~\ref{X}, we examine additional classes of spaces $C_k(X)$ and $C_p(X)$ whose duals admit weak$^\ast$-basic sequences, see Theorem \ref{ZZ}.

\begin{remark}\label{ZZl}
Recall that every metrizable $C_k(X)$ is barrelled; this follows, for example, from \cite[Theorem 10.1.12]{Bonet}.
There also exist non-barrelled spaces $C_k(X)$ for which the conclusion of Theorem~\ref{VVI}(3) holds (see Example~\ref{EU} below).
\end{remark}

\noindent\textbf{(B) Banach-space input and quojections.}
The Josefson--Nissenzweig theorem provides normal sequences in duals of Banach spaces (see \cite{Diestel}).
We shall also use the stronger notion introduced in \cite{S}: a sequence $(y^{*}_{n})\subset E^*$ is \emph{strongly normal} if $\|y_n^*\|=1$ for $n\in \mathbb{N}$ and
\[
\Bigl\{x\in E:\sum_{n}|y^{*}_{n}(x)|<\infty\Bigr\}\ \text{is dense in }E.
\]

On the other hand, Argyros, Dodos and Kanellopoulos proved the following important theorem \cite[Theorem 15]{ADK}:
\begin{theorem}[Argyros--Dodos--Kanellopoulos]\label{II}
The dual space $E^{*}$ of an infinite-dimensional Banach space $E$ admits an infinite-dimensional separable quotient.
\end{theorem}
Combining Theorems~\ref{I} and~\ref{II} yields immediately:
\begin{corollary}\label{III}
The bidual of any infinite-dimensional Banach space has a weak$^\ast$-basic sequence.
\end{corollary}
We extend Corollary~\ref{III} beyond Banach spaces to the class of quojections:
\begin{theorem}\label{VI}
The bidual of any quojection $E$ admits a $w^{*}$-basic sequence.
\end{theorem}
The argument for Theorem~\ref{VI} passes through strict inductive limits and raises, in particular, the question whether the dual of every $(LF)$-space admits a weak$^\ast$-basic sequence.
\begin{problem}\label{IV}
Prove Corollary~\ref{III} without using Theorem~\ref{II}.
\end{problem}
In \cite{K-S-T}, we provided an example of a barrelled space $E$ without an infinite-dimensional separable quotient. By Theorem \ref{I} the weak$^{*}$-dual of $E$ does not admit a weak$^\ast$-basic sequence.
This motivates the following
\begin{problem}\label{Prob:Char}
Characterise those lcs $E$ for which $E^{*}$ admits a weak$^\ast$-basic sequence.
\end{problem}

\begin{remark}
To our knowledge, it is unknown whether the dual $E^{*}$ of any infinite-dimensional separable lcs $E$ contains a weak$^\ast$-basic sequence.
\end{remark}

A Fr\'echet lcs $E$ is called a \emph{quojection} if, for every continuous seminorm $p$ on $E$, the quotient space $E/\ker p$ is normed when endowed with the quotient topology. It is easily seen that $E$ is a quojection if and only if $E = \operatorname{proj}_{n}(E_n, R_n),$
where each $E_n$ is a Banach space and each bonding map $R_n:E_{n+1}\to E_n$ is surjective; in this case we also write $E=\operatorname{quoj}_{n}(E_n,R_n)$.
Any countable product of Banach spaces is a quojection, see \cite[8.4.28]{Bonet} (for example,
$C_{k}(\mathbb{R})\simeq \prod_{n\in\mathbb{N}} C([-n,n])$), but there exist quojections which are not countable products of Banach spaces (Moscatelli, see \cite[Example 8.4.37]{Bonet}).
We refer also to \cite{Moscatelli}, \cite{Vogt}, \cite{Moscatelli-1} and \cite[Examples 8.4.34]{Bonet} for further examples.

Recall finally that $\varphi$ is an $\aleph_{0}$-dimensional vector space with the finest locally convex topology; equivalently, a base at the origin consists of all absolutely convex absorbing sets.
Then every linear functional on $\varphi$ is continuous, and every bounded set $B$ in $\varphi$ is finite-dimensional (that is, $\operatorname{span}(B)$ is finite-dimensional), see \cite[Chapter 2.4]{KKPS} or \cite[0.3.1]{Bonet}.

\section{Proof of Theorem \ref{VVI} and more on spaces $C_k(X)$, $C_p(X)$}\label{X}
By $C_p(X)$ we mean the space of continuous real-valued functions on a Tychonoff space $X$ with the pointwise topology. The space $C_p(X)$ is barrelled if and only if every functionally bounded subset of $X$ is finite, see \cite[Corollary 11.7.6]{Jarchow}.

 Recall that a topological space $X$ is \emph{pseudocompact  } if every real-valued continuous function on $X$ is bounded. It is known that $X$ is not pseudocompact if and only if  $C_p(X)$ contains a (complemented) copy of  the space $\mathbb{R}^{\mathbb{N}}$ if and only if $C_k(X)$ contains a (complemented)  copy of $\mathbb{R}^{\mathbb{N}}$ see  \cite[Section 4]{Ark-2}, or \cite[Theorem 2.6.3]{KKPS}, \cite[Theorem 3.1]{K-S-T}. Then the weak$^\ast$-dual of  $C_p(X)$ contains a copy of the space $(\mathbb{R}^{\mathbb{N}})^{*}=\varphi$; note here that the  Hamel basis of $\varphi$ is also a Schauder basis.

\begin{remark} \label{CD}  In \cite{K-S-T} we provided a  barrelled lcs $E$  without infinite-dimensional separable quotients; hence  $E$ is not isomorphic to any space $C_k(X)$ or $C_p(X)$
by Theorem \ref{VVI}(3), see also \cite[Theorem 2]{Kakol-Saxon-2}.
\end{remark}
We are in the position to prove Theorem \ref{VVI}.  We will need, however, the following result due to \'Sliwa \cite[Theorem 1]{S}.
\begin{theorem}[\'Sliwa]\label{SL}
Any strongly normal sequence $(y_n)$ in the dual $E^{*}$ of a Banach space $E$ contains a weak$^\ast$-basic subsequence.
\end{theorem}
\begin{proof}[Proof of Theorem \ref{VVI}]
(1) By \cite[Corollary 1.10]{KMSZ} the space $C_k(X\times Y)$ has a quotient isomorphic  to $\mathbb{R}^{\mathbb{N}}$ or $c_0$ or $(c_0)_{p}$, where $(c_0)_{p}=\{(x_{n})\in \mathbb{R}^{\mathbb{N}}: x_{n}\rightarrow 0\}$ is endowed with the pointwise topology. Since the above quotient spaces admit Schauder bases, we apply Theorem \ref{I}.

(2) By Warner \cite{warner} the space $C_{k}(X)$ admits a fundamental sequence of bounded sets if and only if for $X$ the following condition holds: Given any sequence $(G_{n})_{n}$ of pairwise disjoint non-empty open subsets of $X$ there is a compact set $K\subset X$ such that $\{n\in\mathbb{N}:K\cap G_{n}\neq\emptyset\}$ is infinite. We say  $X$ is Warner bounded if the above condition holds.

(a) \emph{Assume  that $X$ is Warner bounded.} Then by \cite[Theorem 3.5]{K-S-T} the strong dual of the space $C_{k}(X)$ contains a complemented copy of the space $\ell_{1}$.

(b)\emph{ Assume that  $X$ is not Warner bounded}. Then the strong dual of the space $C_{k}(X)$ contains a complemented copy of $\varphi$ again by \cite[Theorem 3.5]{K-S-T}.

Next, for the both cases we apply Theorem \ref{I} to complete the proof of the item (2).

(3) Recall that $C_k(X)$ is barrelled if and only if  every functionally bounded (=bounding) subset of $X$ is relatively compact, see \cite[Theorem 10.1.20]{Bonet}. If $C_k(X)$  is  a Banach space then $C_{k}(X)$ has a quotient isomorphic to $c_0$ or $\ell_{2}$, see \cite{Rosenthal}.  If $C_k(X)$ is not a Banach space, then $X$ is not pseudocompact. Then $C_k(X)$ contains a complemented copy of $\mathbb{R}^{\mathbb{N}}$. For both cases, Theorem \ref{I} completes the proof.

(4) If $X$ contains an infinite compact subset $K$, the restriction map $T: C_k(X) \rightarrow C(K)$ is continuous and open (\cite[Proposition 2.9]{KMSZ}) and $C(K)$ maps onto $\ell_{2}$ or $c_0$ by a continuous open surjection. Then we apply Theorem \ref{I}. Now assume that all compact subsets of $X$ are finite.  If $X$ is pseudocompact we apply  \cite[Theorem 1]{BKS} to get that $C_k(X)=C_p(X)$ has a quotient with a Schauder basis isomorphic with  $(\ell_{\infty})_p=\{(x_n)\in\mathbb{R}^{\mathbb{N}}: sup_n |x_{n}|<\infty\}$ with the pointwise topology and we again  apply  Theorem \ref{I}. Finally, if $X$ is not pseudocompact, the space $C_k(X)$ contains a complemented copy of $\mathbb{R}^{\mathbb{N}}$,  and Theorem \ref{I} applies.

(5) For every two different points $(x_1, y_1)$ and $(x_2, y_2)$ of $X \times Y$ we have $x_1\neq x_2$ or $y_1\neq y_2$. Thus there exists an element $f\in C(X)$ with $f(x_1)=1$ and $f(x_2)=0$ or $g\in C(Y)$ with $g(y_1)=1$ and $g(y_2)=0$. Hence  $f\otimes {\bf 1}_Y$ or ${\bf 1}_X \otimes g$ separates the points $(x_1, y_1)$ and $(x_2, y_2)$.

Thus the subalgebra $\mathcal {A}:=\mbox{lin}\; \{f\otimes g: (f, g)\in C(X)\times C(Y)\}$ of $C(X \times Y)$ separates points of $X \times Y$. Using the Stone-Weierstrass theorem, we infer that ${\mathcal A}$ is dense in $C(X \times Y)$.

Let $(\mu_{n})$ be a sequence mentioned  in
Theorem \ref{SU} and  let $(\mu_{k_n})$ be a  subsequence of  $(\mu_n)$.  Let $(s_{k_n})$ be a subsequence of $(k_n)$ such that $$\sum_{n=1}^{\infty} \frac{1}{\sqrt{s_{k_n}}}<\infty.$$ Then (applying Theorem \ref{SU}) we derive that for all $(f, g)\in C(X)\times C(Y)$ we have \[ \sum_{n=1}^{\infty} |\mu_{s_{k_n}}(f\otimes g)|\leq \sum_{n=1}^{\infty} \frac{8}{\sqrt{\pi}}\frac{1}{\sqrt{s_{k_n}}}\|f\otimes g\|_{\infty}<\infty.\] It follows that $\sum_{n=1}^{\infty} |\mu_{s_{k_n}}(h)|< \infty$ for every $h\in {\mathcal A}$. Thus the sequence $(\mu_{s_{k_n}})$ is strongly normal, i.e. $\{f\in C(X\times Y): \sum_{n=1}^{\infty} |\mu_{s_{k_n}}(f)|<\infty\}$ is dense in $C(X\times Y)$.

 By  Theorem \ref{SL} every strongly normal sequence $(y^{*}_n)$ in the dual of a Banach space $E$ contains a subsequence that is a weak$^\ast$-basic sequence in the dual of $E$. Using this theorem,  we complete the proof for our case $E=C(X\times Y)$.
 \end{proof}
Corollary \ref{III} applies also to get the following
\begin{corollary}\label{Z}
Let $E$ be an infinite-dimensional lcs containing an infinite-dimensional complemented  Banach subspace $F$. Then the strong bidual $E^{**}_{b}$ of $E$ admits a weak$^\ast$-basic sequence.
\end{corollary}
\begin{proof}
Since $F$ is complemented in $E$, the strong dual $F^{*}_{b}$ of $F$ is complemented in the strong dual $E^{*}_{b}$ of $E$, see \cite[page 287]{Kothe}. Let $P: E^{*}_{b}\rightarrow F^{*}_{b}$ be a continuous linear projection. Applying Theorem \ref{II} for the Banach space $F^{*}_{b}$ and Johnson-Rosenthal \cite[Theorem 2]{JR} we derive that there exists a continuous and open linear surjection $T: E^{*}_{b}\rightarrow G$ onto an infinite-dimensional separable Banach space $G$ with a Schauder basis. Now Theorem \ref{I} applies.
\end{proof}
In \cite[Theorem 32]{Kakol-Saxon} K\c akol and Saxon proved a theorem stating that every infinite-dimensional subspace of the weak$^\ast$-dual of $C_p(X)$ admits an $\aleph_{0}$-dimensional quotient. However, it is still unknown if any space $C_p(X)$ admits an infinite-dimensional separable quotient (even for an infinite compact $X$). Nevertheless, we have the following variant of Theorem \ref{VVI}(3).
\begin{corollary}
If $C_p(X)$ is a barrelled space over infinite Tychonoff $X$, then the dual of $C_{p}(X)$ has a weak$^\ast$-basic sequence.
\end{corollary}
\begin{proof} Consider two cases. (a) $X$ is \emph{pseudocompact}. Then, the compact-open topology on $C(X)$ is dominated by the  Banach topology $\xi$ such that $(C(X),\xi)$ is isomorphic to the Banach space  $C(\beta X)$. Then the identity map
$I: C_p(X)\rightarrow (C(X),\xi)$ has closed graph, so by the closed graph theorem, see \cite[Theorem 4.1.10]{Bonet}, the map $I$ is continuous, implying that $C_p(X)$ is a Banach space. Hence $X$ is finite, a contradiction. (b) $X$ is \emph{not pseudocompact}. Then  $C_{p}(X)$ has a quotient isomorphic with $\mathbb{R}^{\mathbb{N}}$. We use then Theorem \ref{I}.
\end{proof}
Theorem \ref{VVI} may also suggest the following
\begin{problem}\label{WU}
Does $C_p(X)^{*}$ admit a  weak$^\ast$-basic sequence for any infinite compact space $X$?
\end{problem}
\begin{problem}\label{su}
 Does $C_k(X)^{*}$ admit a weak$^\ast$-basic sequence for any infinite Tychonoff space $X$?
\end{problem}
\begin{remark} Note that in Problem \ref{su} it is enough to consider the case when every compact subset of $X$ is finite, and $X$ is pseudocompact. The remaining cases are covered by the previous results.
Therefore, having in mind Theorem \ref{VVI}(4), one may wonder whether for any pseudocompact space $X$ whose every compact subset is finite, there exists in $X$ an infinite discrete subset which is $C^{*}$-embedded. The answer is negative. This example was kindly provided to the first-named author by Professor Reznichenko.
\end{remark}
\begin{example}
In \cite{Reznichenko-2} Tkachenko and Reznichenko quoted the following properties of the space $X$ constructed in \cite{Reznichenko}:
\begin{enumerate}
\item $X$ is a pseudocompact space of cardinality continuum, and each subset of $X$ whose cardinality is less than the continuum is closed and discrete.
\item  $\beta X$ is homeomorphic to a Tychonoff cube of uncountable weight and is a Dugundji compactum. 
\item The closure of every countable subset of X in $\beta X$ is a metrizable compact space.
\end{enumerate}  Hence  $X$  contains no countable discrete subsets $C^{*}$-embedded in $X$ and every compact subset of $X$ is finite.
\end{example}
Note also that the conclusion of Theorem \ref{VVI} (3) holds also for nonseparable spaces $C_k(X)$ not being barrelled.
\begin{example}\label{EU}
Let $X$ be a pseudocompact space such that $X$ contains a copy of $\mathbb{N}$ and is contained in $\beta\mathbb{N}$, and each compact subset of $X$ is finite, see \cite{Haydon}. Then $C_k(X)=C_p(X)$ is not barrelled. Indeed,  if $X$ is pseudocompact, then  $C_k(X)$ admits a stronger Banach topology generated by the norm $\|f\|= sup_{x\in X}|f(x)|.$  The closed  graph theorem  between barrelled and Banach spaces \cite[Theorem 10.1.10]{Bonet} implies that $X$ must be finite, a contradiction.  Since $X$ is pseudocomopact containing an infinite discrete subset $C^{*}$-embedded into $X$, we apply again   \cite[Theorem 1]{BKS} to deduce  that $C_p(X)$ has a quotient isomorphic to  $(\ell_{\infty})_p$ with a Schauder basis, see also \cite[Proposition 3.1 (4)]{KMSZ} for that case. Now Theorem \ref{I} applies. Note also that $C_k(X)$ is not separable, otherwise $X$ would admit a weaker metrizable topology, and then by \cite[IV.5.10, page 156]{Ark} we conclude that $X$ is metrizable and compact, what is impossible.
\end{example}
Applying previous facts, we  note the following
\begin{corollary}
Let $X$ be an infinite Tychonoff space. If $C_k(X)$ does not admit an infinite-dimensional quotient with a Schauder basis, then $X$ is pseudocompact, and every compact subset of $X$ is finite. Example \ref{EU} shows that the converse implication fails in general.
\end{corollary}
\begin{corollary}\label{true}
Assume that for an infinite Tychonoff space $X$, the space $C_p(X)$ is separable. Then
\begin{enumerate}
\item $C_k(X)$ admits an infinite-dimensional quotient isomorphic to $\mathbb{R}^{\mathbb{N}}$ or $c_0$ or $\ell_{2}.$
\item $C_p(X)$ admits an infinite-dimensional quotient isomorphic to $\mathbb{R}^{\mathbb{N}}$ or $(c_0)_p=\{(x_n)\in\mathbb{R}^{\mathbb{N}}: x_{n}\rightarrow 0\}$ endowed with the product topology from $\mathbb{R}^{\mathbb{N}}$.
\end{enumerate}
\end{corollary}
\begin{proof}
Since $C_p(X)$ is separable, $X$ admits a weaker metrizable separable topology, \cite[Corollary 4.2.2]{mccoy}.

(1) \emph{If $X$ is pseudocompact}, then  applying again \cite[IV.5.10, page 156]{Ark} we deduce that $X$ is a metrizable compact space, so $C_{k}(X)$ has a quotient isomorphic to  $c_0$.  \emph{If $X$ is not pseudocompact}, then  $C_{k}(X)$ contains a complemented copy of $\mathbb{R}^{\mathbb{N}}$.

(2) If $X$ is \emph{not pseudocompact}, $C_p(X)$ contains a complemented copy of $\mathbb{R}^{\mathbb{N}}$.  If $X$ \emph{is pseudocompact}, $X$ is compact metrizable. Then $C_p(X)$ contains a complemented copy of $(c_0)_p$, see \cite{BKS1}.
\end{proof}
\begin{problem}\label{RR}
Assume that $C_k(X)$ ($C_{p}(X)$, respectively) has an infinite-dimensional separable quotient. Is it true that then $C_{k}(X)$ ($C_{p}(X)$, respectively) admits an infinite-dimensional quotient with a Schauder basis?
\end{problem}
In \cite{BKS} we investigated the problem when the space $C_p(X)$ has an infinitely dimensional metrizable quotient.
 This also motivates  Theorem \ref{ZZ} related to Problem \ref{WU}, Problem \ref{RR}  and Corollary \ref{true}.
\begin{theorem}\label{ZZ}
Let $X$ be an infinite Tychonoff space.
\begin{enumerate}
\item If the space  $C_p(X)$  has an infinite-dimensional  metrizable quotient $C_p(X)/Y$, then $C_p(X)/Y$ admits a Schauder basis and  the dual  $C_p(X)^{*}$ has a weak$^\ast$-basic sequence.
\item If the space $C_k(X)$ has an infinite-dimensional metrizable quotient $C_k(X)/Y$, then $C_k(X)$ admits also an infinite-dimensional metrizable quotient $C_k(X)/Z$ with a Schauder basis and    $C_k(X)^{*}$ has a weak$^\ast$-basic sequence.
\end{enumerate}
\end{theorem}
For the proof, we need the following two additional lemmas, which might be known but hard to localize.
\begin{lemma}\label{GG}
A linear subspace $Y$  of  $\mathbb{R}^{\mathbb{N}}$ is dense in  $\mathbb{R}^{\mathbb{N}}$ if and only if for every finite set $F\subset \mathbb{N}$ the projection $\pi_{F}: Y\rightarrow\mathbb{R}^{F}$ is surjective, where $\pi_{F}: \mathbb{R}^{\mathbb{N}}\rightarrow\mathbb{R}^{F}$ is the canonical projection given by $\pi_{F}(x)=(x_{n})_{n\in F}$.
\end{lemma}
\begin{proof}
Let $Y\subset\mathbb{R}^{\mathbb{N}}$  be dense. Since $\pi_{F}$ is continuous,  $\pi_{F}(Y)$ is  dense and finite-dimensional in $\pi_{F}(\mathbb{R}^{\mathbb{N}})=\mathbb{R}^{F}$. Hence $\pi_{F}$ is a surjection. Conversely, assume  $\pi_{F}(Y)=\mathbb{R}^{F}$ for every finite $F\subset\mathbb{N}$. Let $U=\pi^{-1}_{F}(V)$ be an open basic set in the product topology of $\mathbb{R}^{\mathbb{N}}$ with open $V\subset\mathbb{R}^{F}$ and some finite  $F\subset\mathbb{N}$. Then
$V\cap Y\neq\emptyset,$  so $Y$ is dense.
\end{proof}
\begin{lemma}\label{SS}
Every dense vector subspace $Y$ of $\mathbb{R}^{\mathbb{N}}$ admits a Schauder basis.
\end{lemma}
\begin{proof}
We construct by induction  a sequence $(b_{n})$ in $Y$ with the condition $\pi_{k}(b_n)=\delta_{kn}$ for all $1\leq k \leq n$.  For $n=1$ by Lemma \ref{GG} we have $\pi_{1}(Y)=\mathbb{R}$, so we find $b_1\in Y$ with $\pi_{1}(b_1)=1.$ Assume $b_1,b_2,\dots b_{n-1}$ have been found. By Lemma \ref{GG} the map $\pi_{F}: Y\rightarrow\mathbb{R}^{n}$ is surjective for $F=\{1,2,\dots, n\}$.  Hence there exists $b_n\in Y$ such that
$$\pi_{F}(b_n)=(0,0,\dots, 1).$$ We show that $(b_n)_n$ is a Schauder basis. Let $y=(y_n)\in Y$. Define  the scalars $(a_{n})_{n}$ as follows:
$$a_1=y_1, \,\, a_n=y_n-\sum_{k=1}^{n-1}a_k\pi_{n}(b_{k}),\,\, n>1.$$
Set $S_{N}=\sum_{n=1}^{N}a_nb_n.$ Fix $m\in\mathbb{N}$. For all $N\geq m$ we have
$$ \pi_{m}(S_{N})=\sum_{n=1}^{N}a_n\pi_{m}(b_n)=\sum_{n=1}^{m}a_{n}\pi_{m}(b_n),$$
since $\pi_{m}(b_{n})=0$ for all $n> m$. Note that $\pi_{m}(S_N)=y_m$ for all $N\geq m$. Indeed, for $m=1$ we have $\pi_{1}(S_N)=a_1\pi_{1}(b_1)=a_1=y_1.$ Assume that the claim holds for all $k<m.$ Then
$$\pi_{m}(S_N)=\sum_{n=1}^{m-1}a_n\pi_{m}(b_n)+a_m\pi_{m}(b_m)$$
$$=\sum_{n=1}^{m-1}a_n\pi_{m}(b_n)+ (y_m-\sum_{k=1}^{m-1}a_k\pi_{m}(b_k)).$$
This implies that $\pi_{m}(S_N)\rightarrow_N y_m$ for every $m\in\mathbb{N}$. Hence $S_N\rightarrow_N y$ in the product topology. The uniqueness of the coefficients follows from the triangular structure, $a_{n}$ is uniquely determined by $\pi_{n}(y)$ and the preceding coefficients.

Finally we show (by induction on $n$) that the basic projections $b_{n}^{*}: Y\rightarrow\mathbb{R}$ defined by $b_{n}^{*}(y)=a_n$ are continuous. For $n=1, \, b_{1}^{*}(y)=\pi_{1}(y).$ Since $\pi_{1}$ is continuous, its restriction to $Y$ is continuous. For $n> 1$
$$ b_{n}^{*}(y)=\pi_{n}(y)-\sum_{k=1}^{n-1}b_{k}^{*}(y)\pi_{n}(b_k),$$
which shows the continuity of $b_{n}^{*}$.
\end{proof}
Theorem  \ref{ZZ} combined with Theorem \ref{I} motivate the next
\begin{problem}
Is it true that  $C_p(X)$ with an infinite-dimensional separable quotient admits also an infinite-dimensional metrizable quotient?
\end{problem}
The answer is positive if $C_p(X)$ itself is separable, see Corollary \ref{true}.
\begin{proof}[Proof of Theorem \ref{ZZ}]
(1) Let $C_p(X)/Y$ be  an infinite-dimensional metrizable quotient of $C_p(X)$. Since $C_p(X)$ and then also  $C_p(X)/Y$ carry its weak topologies, respectively, the metrizable space $C_p(X)/Y$ embeds into $\mathbb{R}^{\mathbb{N}}$. Consequently, $C_p(X)/Y$ is separable.  We may assume that $E=C_p(X)/Y$ is a dense subspace of $\mathbb{R}^{\mathbb{N}}$ (since every closed infinite-dimensional vector subspace  of $\mathbb{R}^{\mathbb{N}}$ is isomorphic to $\mathbb{R}^{\mathbb{N}}$, see
\cite[Corollary 2.6.5]{Bonet}).
Applying Lemma \ref{SS} we conclude that $C_p(X)/Y$ has a Schauder basis. Finally, we apply Theorem \ref{I}.

(2)  If every compact subset of $X$ is finite, then $C_p(X)=C_k(X)$ and item  (1) applies. If $X$ contains an infinite compact subset, we know already that $C_k(X)$ has a quotient isomorphic to $c_0$ or $\ell_{2}$, and then again Theorem \ref{I} applies.
\end{proof}
We apply the above arguments to derive the following

\begin{corollary}
 If $X$ is a Tychonoff space such that the dual of $C_p(X)$ admits a weak$^{*}$-basic sequence, then the same property holds for the dual of $C_k(X)$.
\end{corollary}
Being motivated by the above results related with spaces $C_k(X)$ and $C_p(X)$ we may pose the following
\begin{problem}
Does an infinite-dimensional metrizable and separable lcs admit an infinite-dimensional quotient with a Schauder basis?
\end{problem}
It is known that a Tychonoff space $X$ is separable if and only if $C_p(X)$ admits a weaker metrizable (separable) locally convex topology, say $\xi$, see \cite[Corollary 4.3.3]{mccoy}. By \cite[Theorem 1.5]{Kakol-Leiderman} the completion of $(C_p(X),\xi)$ is isomorphic to $\mathbb{R}^{\mathbb{N}}$. Applying Lemma \ref{SS} we note the following
\begin{corollary}\label{WW}
A Tychonoff space $X$ is separable if and only if  $C_p(X)$  admits a weaker metrizable locally convex topology with a Schauder basis.
\end{corollary}

\section{Proof of Theorem \ref{VI}} \label{XX}
Let $E$ be a vector space, and let $(E_{n},\tau_{n})$ be a
strictly increasing sequence of vector subspaces of $E$ covering $E$, each
$E_{n}$ is a Fr\'echet (Banach) space  endowed with the topology $\tau_{n}$, such that
$\tau_{n+1}|E_{n}\leq\tau_{n}$ for each $n\in\mathbb{N}$. Then on
$E$ there exists the finest locally convex topology $\tau$ such
that $\tau|E_{n}\leq\tau_{n}$ for each $n\in\mathbb{N}$, and then
the space  $(E,\tau)$ is called  an  \textit{$(LF)$-space}
 (an \textit{$(LB)$-space}).

 If $\tau_{n+1}|E_{n}=\tau_{n}$ for each $n\in\mathbb{N}$, then
$\tau|E_{n}=\tau_{n}$ for each $n\in\mathbb{N}$, and $(E,\tau)$ is called the \textit{strict inductive limit} of $(E_{n},\tau_{n})$, see for example \cite{Bonet}. Recall that every $(LF)$-space is barrelled, see \cite[Proposition 4.2.6]{Bonet}. The classical spaces $D^{m}(\Omega)$ (strict $(LB)$-space),  $D(\Omega)$, $D(\mathbb{R}^{N})$  (strict $(LF)$-spaces) illustrate the above classes.

In order to prove Theorem \ref{VI} we  need also the following proposition. 
\begin{proposition} \label{VII}  The weak$^\ast$-dual $E^{*}$ of a strict  $(LF)$-space contains  a basic  sequence. In fact, it contains a complemented copy of  $\mathbb{R}^{\mathbb{N}}$.
\end{proposition}
\begin{proof} By \cite[Theorem 1]{Saxon} the strict $(LF)$-space $E$ contains a complemented
copy of the space $\varphi$. Hence there exist continuous linear maps
$i:\varphi\rightarrow E$ and $P:E\rightarrow\varphi$ such that $P\circ i=id_{\varphi}.$

Passing to adjoints, we obtain continuous maps
$P^*:(\varphi^*,\sigma(\varphi^*,\varphi))\rightarrow (E^*,\sigma(E^*,E))$ and
$i^*:(E^*,\sigma(E^*,E))\rightarrow(\varphi^*,\sigma(\varphi^*,\varphi))$
satisfying $i^*\circ P^*=id_{\varphi^*}.$ Hence $P^*$ is a topological embedding and $P^*(\varphi^*)$ is complemented (in particular, closed) in the weak$^\ast$-dual of $E$.

Finally, $\varphi^*$ is isomorphic to $\mathbb{R}^{\mathbb{N}}$ and $\sigma(\varphi^*,\varphi)$ corresponds to the product topology on $\mathbb{R}^{\mathbb{N}}$. The unit vectors form a Schauder basis
of $\mathbb{R}^{\mathbb{N}}$ in this topology, hence their image under $P^*$ is a $w^*$-basic sequence in $E^*$. This also shows that the weak$^\ast$-dual of $E^*$ contains a complemented copy of $\mathbb{R}^{\mathbb{N}}$.
\end{proof}
In \cite[Theorem 3]{Saxon-1}, Saxon and Narayanaswami proved that every $(LF)$-space $E$  admits an infinite-dimensional separable quotient. If $E$ is a strict $(LF)$-space, then $E$ has a complemented copy of $\varphi$ (\cite[Theorem 1]{Saxon}). The above facts, combined with Theorem \ref{I} motivate
\begin{problem}
Does the dual of every  $(LF)$-space admit a weak$^\ast$-basic sequence?
\end{problem}
\begin{proof}[Proof of Theorem \ref{VI}]
If $E$ is a Banach space, we apply Corollary \ref{III}.

Now assume that $E$ is a Fr\'echet lcs which is a quojection and $E$ is not normable. Since $E$ is a quojection, the space  $E$ is the strict projective limit of Banach spaces
$E= quoj_n(E_n, R_n).$  The dual $E^{*}$ equipped with the strong topology $\beta(E^{*},E)$ is then the strict inductive  limit of the sequence $(E_{n}^{*},\beta(E_{n}^{*},E_{n}))$ of Banach spaces by
\cite[Proposition 8.4.30]{Bonet}.

Note that the  sequence $(E_{n}^{*})$ is indeed \emph{strictly increasing}. In fact,  since the topology $\beta(E^{*},E)$ is barrelled (as the inductive limit topology of Banach spaces $(E_{n}^{*},\beta(E_{n}^{*},E_{n}))$, see  \cite[Proposition 4.2.6]{Bonet}) and non-normable (since $E$ is not a Banach space),  we can select a strictly increasing subsequence of $(E_{n}^{*})$, covering $E^{*}$. Otherwise, we could cover $\beta(E^{*},E)$ by a sequence of closed absolutely convex sets of the form $(nS)$, where $S$ is a bounded set in $\beta(E^{*},E)$. Hence $S$ would be a barrel in $\beta(E^{*},E)$, and then by the barreledness of $\beta(E^{*},E)$, this space would have a bounded neighbourhood of zero yielding the normability of $\beta(E^{*},E)$, a contradiction. Hence $(E^{*}, \beta(E^{*},E))$ is a strict $(LB)$-space. Finally, Proposition \ref{VII} completes the proof.
\end{proof}

\backmatter





\section*{Declarations}


\begin{itemize}
\item Funding: Not applicable.
\item Conflict of interest/Competing interests: The authors have no competing interests.
\item Ethics approval and consent to participate: Not applicable.
\item Consent for publication: All authors have read and approved the final manuscript and consent to its submission and publication.
\item Data availability: Not applicable.
\item Materials availability: Not applicable.
\item Code availability: Not applicable.
\item Author contribution: All authors contributed to the conception of the work. All authors wrote and approved the final manuscript.
\end{itemize}








\bibliography{sn-bibliography}

\end{document}